\documentclass[12pt]{amsart}

\advance\textheight\topskip
\allowdisplaybreaks

\usepackage{mathabx}
\newcommand{\Smash}{\mathbin{\hash}}
\newcommand{\Smashfn}{\Smash}

\newcommand{\Braid}{\mathbin{\mbox{\large${\Join}$}}}

\usepackage{amsfonts,amssymb}
\usepackage{times}
\usepackage[mathscr]{eucal}
\usepackage{mathptmx}
\usepackage{amsthm}

\newcommand{\DDelta}{\Delta}
\newcommand{\tensorA}{\mathbin{{\otimes}_{_{\Mat_p(\oC)\!\!}}}}
\newcommand{\E}{\mathsf{e}}
\newcommand{\one}{\pmb{1}}

\newcommand{\XX}{\boldsymbol{\mathscr{A}}}

\newcommand{\pS}{\mathsf{s}}

\newcommand{\Matrixi}[1]{\mbox{\small$\begin{pmatrix}#1\end{pmatrix}$}}
\newcommand{\Matrixii}[1]{\mbox{\footnotesize$\begin{pmatrix}#1\end{pmatrix}$}}

\newcommand{\Dz}{\dd}

\newcommand{\zM}{{Z}}
\newcommand{\DzM}{{D}}

\newcommand{\up}[1]{^{{\scriptscriptstyle(#1)}}}

\newcommand{\zero}{_{_{(0)}}}
\newcommand{\mone}{_{_{(-1)}}}

\newcommand{\adjoint}{\kern2pt\mbox{\footnotesize$\blacktriangleright$}\kern2pt}
\newcommand{\leftreg}{\kern1pt{\rightharpoonup}\kern1pt}
\newcommand{\rightreg}{\kern1pt{\leftharpoonup}\kern1pt}

\newcommand{\leftregchi}{\kern1pt{\rightharpoonup}_{\!\!\!\!_{\chi}}\kern1pt}

\newcommand{\doteta}{\kern2pt{\cdot}_{\!_{\eta}}\kern1pt}

\newcommand{\Czd}{\oC_{\q}[z,\Dz]}
\newcommand{\Mat}{\mathrm{Mat}}
\renewcommand{\kappa}{\varkappa}
\newcommand{\qfac}[1]{[#1]!\,}
\newcommand{\qint}[1]{[#1]}

\newcommand{\leftii}{\mathbin{\mbox{\footnotesize${\vartriangleright}$}}}

\newcommand{\Dp}[1]{\,_{\phantom{h}}^{\underline{#1\kern-.5pt}\kern.5pt}}

\newcommand{\Sinv}{{S^*}^{-1}}

\newcommand{\Star}{\mathbin{{\star}}}

\usepackage[dvips, bookmarks=false, hyperfigures=false,
setpagesize=false]{hyperref}

\newcommand{\id}{\mathrm{id}}
\newcommand{\eval}[2]{\langle#1,\,#2\rangle\,}

\newcommand{\DD}{\mathscr{D}}
\newcommand{\HD}{\mathscr{H}}
\newcommand{\HH}{\boldsymbol{\mathscr{H}}}
\newcommand{\HHH}{\pmb{\textsf{H}}}
\newcommand{\bDDB}{\kern4pt\overline{\kern-3pt\mathscr{D}(B)\kern-3pt}\kern4pt}

\newcommand{\bref}[1]{\textbf{\ref{#1}}}

\renewcommand{\geq}{\geqslant}
\renewcommand{\leq}{\leqslant}

\newcommand{\tensor}{\otimes}

\newcommand{\q}{\mathfrak{q}}

\newcommand{\UresSL}[1]{\overline{\mathscr{U}}_{\q} s\ell(#1)}
\newcommand{\HresSL}[1]{\overline{\mathscr{H}}_{\q} s\ell(#1)}

\newcommand{\mfrac}[2]{\raisebox{.8pt}{\mbox{\small$\displaystyle\frac{#1}{#2}$}}}
\newcommand{\ffrac}[2]{\raisebox{.5pt}{\mbox{\footnotesize$\displaystyle\frac{#1}{#2}$}}}
\newcommand{\fffrac}[2]{\raisebox{.9pt}{\mbox{\scriptsize$\displaystyle\frac{#1}{#2}$}}}
\newcommand{\half}{%
  \mathchoice{\ffrac{1}{2}}{\frac{1}{2}}{\frac{1}{2}}{\frac{1}{2}}}

\newcommand{\qbin}[2]{\mathchoice%
  {\qbinm{#1}{#2}}{\qbinmm{#1}{#2}}%
  {\qbinmm{#1}{#2}}{\qbinmm{#1}{#2}}}
\newcommand{\qbinm}[2]{\mbox{\footnotesize$\displaystyle
    \genfrac{[}{]}{0pt}{}{#1}{#2}$}}
\newcommand{\qbinmm}[2]{\genfrac{[}{]}{0pt}{}{#1}{#2}}

\newcommand{\dd}{\partial}

\newcommand{\oC}{\mathbb{C}}

\newcommand{\oZ}{\mathbb{Z}}

\numberwithin{equation}{section}
\makeatletter
\@addtoreset{equation}{section}
\@addtoreset{subsubsection}{section}

\def\@secnumfont{\bfseries}
\def\subsubsection{\@startsection{subsubsection}{3}%
  \z@{.5\linespacing\@plus.7\linespacing}{-.5em}%
  {\normalfont\bfseries}}
\def\paragraph{\@startsection{paragraph}{4}%
  \z@\z@{-\fontdimen2\font}%
  \normalfont\bfseries}
\def\subparagraph{\@startsection{subparagraph}{5}%
  \z@\z@{-\fontdimen2\font}%
  \normalfont\bfseries}

\makeatother

\swapnumbers
\newtheorem{Thm}[subsection]{Theorem}
\newtheorem{thm}[subsubsection]{Theorem}

\newtheorem{lemma}[subsubsection]{Lemma}

\theoremstyle{definition}

\begin{document}
\vspace*{-\baselineskip}

\title[Yetter--Drinfeld module algebras]{Heisenberg double
  $\HD(B^*)$ as a braided commutative Yetter--Drinfeld module algebra
  over the Drinfeld double}

\author[Semikhatov]{A.M.~Semikhatov}

\address{Lebedev Physics Institute
  \hfill\mbox{}\linebreak \texttt{ams@sci.lebedev.ru}}

\begin{abstract}
  We study the Yetter--Drinfeld $\DD(B)$-module algebra structure on
  the Heisenberg double $\HD(B^*)$ endowed with a ``heterotic'' action
  of the Drinfeld double $\DD(B)$.  This action can be interpreted in
  the spirit of Lu's description of $\HD(B^*)$ as a twist of~$\DD(B)$.
  In terms of the braiding of Yetter--Drinfeld modules, $\HD(B^*)$ is
  braided commutative.  By the Brzezi\'nski--Militaru theorem,
  $\HD(B^*)\Smash\DD(B)$ is then a Hopf algebroid over~$\HD(B^*)$. \
  For $B$ a particular Taft Hopf algebra at a $2p$th root of unity,
  the construction is adapted to yield Yetter--Drin\-feld module
  algebras over the $2p^3$-dimensional quantum group~$\UresSL2$.  In
  particular, it follows that $\Mat_p(\oC)$ is a braided commutative
  Yetter--Drin\-feld $\UresSL2$-module algebra and $\Mat_p(\UresSL2)$
  is a Hopf algebroid over $\Mat_p(\oC)$.
\end{abstract}

\maketitle
\thispagestyle{empty}

\section{Introduction}
For a Hopf algebra $B$, the Heisenberg double $\HD(B^*)$ is the smash
product $B^*\Smash B$ with respect to the left regular action
$b\leftreg \beta=\eval{\beta''}{b}\beta'$ of $B$ on $B^*$; the
composition in $\HD(B^*)$ is given~by
\begin{equation}\label{H-comp}
  (\alpha\Smash a)(\beta\Smash b)=
  \alpha(a'\leftreg\beta)\Smash a'' b,
  \qquad
  \alpha,\beta\in B^*,\quad
  a,b\in B.
\end{equation}
Let $\DD(B)$ be the Drin\-feld double of $B$, with its elements
written as $\mu\tensor m$, where $\mu\in B^*$ and $m\in B$; the
composition in $\DD(B)$ is $(\mu\tensor m)(\nu\tensor n) =
\mu(m'\leftreg \nu\rightreg S^{-1}(m''')) \tensor m''n$ (and the
co\-algebra structure is that of $B^{*\mathrm{cop}}\tensor B$).  We
define a $\DD(B)$ action on~$\HD(B^*)$ as
\begin{align}\label{the-action}
  (\mu\tensor m)\leftii(\beta\Smash b)
  =\mu'''(m'\leftreg\beta)\Sinv(\mu'')
  \Smash\bigl((m'' b S(m'''))\rightreg\Sinv(\mu')\bigr),\\[-2pt]
  \notag
  \mu\tensor m\in\DD(B),\quad \beta\Smash b\in\HD(B^*),
\end{align}
where $b\rightreg\mu=\eval{\mu}{b'}b''$ is the right regular action of
$B^*$ on $B$ (and $\eval{\;}{\,}$ is the evaluation).

\begin{Thm}\label{thm:YD}
  For a Hopf algebra $B$ with bijective antipode, $\HD(B^*)$ endowed
  with action~\eqref{the-action} and the coaction
  \begin{equation}\label{delta}
    \delta:
    \begin{aligned}
      \HD(B^*)&\to\DD(B)\tensor\HD(B^*)\\[-2pt]
      \beta\Smash b&\mapsto(\beta''\tensor b')
      \tensor (\beta'\Smash b'')
    \end{aligned}
  \end{equation}
  is a \textup{(}left--left\textup{)} Yetter--Drin\-feld
  $\DD(B)$-module algebra.
\end{Thm}

By a Yetter--Drin\-feld module algebra we mean a module comodule
algebra that is also a Yetter--Drin\-feld module, i.e., a
compatibility condition between the action and the coaction holds in
the form\pagebreak[3]
\begin{equation}\label{eq:YD}
  (M'\leftii A)\mone M''\tensor(M'\leftii A)\zero
  =M' A\mone\tensor(M''\leftii A\zero),
\end{equation}
where, in our case, $M\in\DD(B)$ and $A\in\HD(B^*)$.\footnote{For a
  Hopf algebra $H$ and a left $H$-comodule $X$, we write the coaction
  $\delta:X\to H\tensor X$ as $\delta (x)=x\mone\!\tensor x\zero$;
  then the comodule axioms are $\eval{\varepsilon}{\!x\mone}x\zero =x$
  and $x\mone'\!\tensor x\mone''\!\tensor x\zero =x\mone\!\tensor
  x\zero{}\mone\!\tensor x\zero{}\zero$.}

We recall from~\cite{[CFM],[BM],[CGW]} that for a Hopf algebra $H$, a
left $H$-module and left $H$-comodule algebra $X$ is said to be
\textit{braided commutative} (or $H$-commutative) if
\begin{equation}\label{braided-comm}
  y x = (y\mone\leftii x) y\zero,\qquad x,y\in X.
\end{equation}
Also, for any two (left--left) Yetter--Drin\-feld $H$-module algebras
$X$ and $Y$, their \textit{braided product} $X\Braid Y$ is defined as
the tensor product with the composition
\begin{equation}\label{braided-prod}
  (x\Braid y)(v\Braid u)
  = x(y\mone\leftii v)\Braid y\zero u,
  \quad x,v\in X,\quad y,u\in Y.
\end{equation}
(This gives a Yetter--Drin\-feld module algebra.)

\begin{Thm}\label{thm:braided-comm}
  $\HD(B^*)$ is a braided \textup{(}$\DD(B)$-\textup{)} commutative
  algebra.  Moreover, $\HD(B^*)$ is the braided product
  \begin{equation*}
    \HD(B^*) = B^{*\mathrm{cop}}\Braid B,
  \end{equation*}
  where $B^{*\mathrm{cop}}$ and $B$ are \textup{(}braided
  commutative\textup{)} Yetter--Drin\-feld $\DD(B)$ module algebras by
  restriction, i.e., with the $\DD(B)$ action
  \begin{equation*}
    (\mu\tensor m)\leftii\beta=\mu''(m\leftreg\beta)\Sinv(\mu'),
    \qquad
    (\mu\tensor m)\leftii b
    =(m' b S(m''))\rightreg\Sinv(\mu)
  \end{equation*}
  and coaction $\delta:\beta\mapsto(\beta''\tensor 1)\tensor\beta'$, \
  $\delta: b\mapsto (\varepsilon\tensor b')\tensor b''$ \
  \textup{(}$\beta\in B^*$, $b\in B$\textup{)}.
\end{Thm}

\subsubsection{}\label{BM-1st}As a corollary, the
Brzezi\'nski--Militaru theorem~\cite{[BM]} then ``provides one with a
rich source of examples of bialgebroids.''  In particular, \textit{for
  any Hopf algebra $B$ with bijective antipode, the ``quadruple''
  $\HD(B^*)\Smash\DD(B)$}, where the smash product is defined with
respect to action~\eqref{the-action}, \textit{is a Hopf algebroid over
  $\HD(B^*)$}.

\subsubsection{A ``pseudoadjoint'' interpretation
  of~\eqref{the-action}}The $\DD(B)$-action~\eqref{the-action} first
appeared in~\cite{[S-H]}.  To borrow a popular term from string
theory~\cite{[GHMR]} (where it was also a borrowing originally), this
action may be termed ``heterotic'' because it is constructed by
combining left and right $\DD(B)$ actions, as we describe
in~\bref{both} (and the heterotic string famously combines ``left''
and ``right'').  Or because~\eqref{the-action} ``cross-breeds''
regular and adjoint actions.

Trying to quantify how ``far'' \eqref{the-action} is from the adjoint
action, we arrive at a useful interpretation of our ``heterotic''
action by extending Lu's description of the product on $\HD(B^*)$ as a
twist of the product on~$\DD(B)$~\cite{[Lu-double]}.  The two
algebraic structures, $\DD(B)$ and $\HD(B^*)$, are defined on the same
vector space $B^*\tensor B$, and the product~\eqref{H-comp} in
$\HD(B^*)$, temporarily denoted by $\Star$, can be written as
\begin{equation}\label{Star-eta}
  M\Star N=M' N' \eta(M'',N''),\qquad M,N\in\DD(B)
\end{equation}
for a certain $2$-cocycle $\eta:\DD(B)\tensor\DD(B)\to k$
\cite{[Lu-double]}.\pagebreak[3] In the same vein, the $\DD(B)$ action
on $\HD(B^*)$ in~\eqref{the-action} can be rewritten in the
``pseudoadjoint'' form
\begin{equation}\label{pseudo-ad-2}
  (M,A)\mapsto M'\Star A\Star\pS(M''),\qquad M\in\DD(B),\quad
  A\in\HD(B^*),
\end{equation}
where $\pS(M)=\eta(M',M'')S(M''')$.  Some ``antipode-like'' properties
of $\pS$ allow independently verifying that the right-hand side here
\textit{is} an action, as we show in~\bref{thm:star}, where further
details are given.

\subsubsection{}\label{Heis-n}
The Heisenberg double $\HD(B^*) = B^{*\mathrm{cop}}\Braid B$ can be
regarded as the lowest term, $\HD(B^*)=\HH_{2}$, in a series of
\textit{Heisenberg $n$-tuples}, or \textit{chains}
$\HH_n$\,---\,the Yetter--Drin\-feld $\DD(B)$-modules
\begin{align*}
  \HH_{2n} &=
  B^{*\mathrm{cop}}\Braid B\Braid B^{*\mathrm{cop}}\Braid B
  \Braid\dots\Braid B,\\
  \HH_{2n+1} &=
  B^{*\mathrm{cop}}\Braid B\Braid B^{*\mathrm{cop}}\Braid B
  \Braid \dots \Braid B\Braid B^{*\mathrm{cop}}
\end{align*}
(with $2n$ and $2n+1$ factors),
with the relations
\begin{align}\label{nn-braiding}
  b[2i]\,\beta[2j+1] &= (b'\leftreg\beta)[2j+1]\,b''[2i]
  \quad\text{for all $i$ and $j$},
  \\
  \intertext{(where $B^{*\mathrm{cop}}\to B^{*\mathrm{cop}}[2j+1]$ and
    $B\to B[2i]$ are the morphisms onto the respective factors,
    and we omit ${}\Braid{}$ for simplicity), and}
  \alpha[2i+1]\,\beta[2j+1]&=(\alpha'''\beta\Sinv(\alpha''))[2j+1]\,
  \alpha'[2i+1],\quad i\geq j\\[-2pt]
  a[2i]\, b[2j]&=(a' b S(a''))[2j]\, a'''[2i], \quad i\geq j,
\end{align}
where $a,b\in B$, $\alpha,\beta\in B^{*\mathrm{cop}}$.

\subsection{}
As regards the popular subject of Yetter--Drin\-feld modules, we note
Refs.~\cite{[Y],[LR],[RT],[Mont],[Sch],[CFM]}.  Heisenberg
doubles~\cite{[AF],[RSts],[Sts],[Lu-double]}, among various smash
products, have attracted some attention, notably in relation to Hopf
algebroid constructions~\cite{[Lu-alg],[P],[BM]} (the basic
observation being that $\HD(B^*)$ is a Hopf algebroid over
$B^*$~\cite{[Lu-alg]}) and also from various other standpoints and for
different purposes~\cite{[K],[Ka],[Mi],[IP]}.  (A relatively recent
paper where Yetter--Drin\-feld-like structures are studied in relation
to ``smash'' products is~\cite{[BB]}.)

\subsection{}The above results are proved quite straightforwardly.
The proofs are given in Sec.~\ref{sec:proofs}; there, $B$ denotes a
Hopf algebra with bijective antipode.  When we pass to an example in
Sec.~\ref{sec:sl2}, $B$ becomes a particular Taft Hopf algebra.

\subsection{}\label{sec:Uq}The example worked out in
Sec.~\ref{sec:sl2} is that of the $2p^3$-dimensional quantum group
$\UresSL2$ at the $2p$th root of unity
\begin{equation*}
  \q=e^{\frac{i\pi}{p}},
\end{equation*}
($p=2,3,\dots$).  This is the Hopf algebra with generators $E$, $K$,
and $F$ and the relations
\begin{gather*}
  KEK^{-1}=\q^2E,\quad
  KFK^{-1}=\q^{-2}F,\quad
  [E,F]=\ffrac{K-K^{-1}}{\q-\q^{-1}},
  \\
  E^{p}=F^{p}=0,\quad K^{2p}=1,
\end{gather*}
and the Hopf algebra structure $\Delta(E)= E\tensor K + 1\tensor E$,
$\Delta(K)=K\tensor K$, $\Delta(F)=F\tensor 1 + K^{-1}\tensor F$,
$\varepsilon(E)=\varepsilon(F)=0$, $\varepsilon(K)=1$, $S(E)=-E K^{-1}$,
$S(K)=K^{-1}$, $S(F)=-K F$.

\subsubsection{}$\UresSL2$ is ``almost'' the Drin\-feld double of a
$4p^2$-dimen\-sional Taft Hopf algebra $B$, more precisely, a
``truncation'' of the double obtained by taking a quotient and then
restricting to a subalgebra.  This close kinship of $\UresSL2$ to a
Drinfeld double extends to the ``Heisenberg side'': it turns out that
the pair $(\DD(B),\HD(B^*))$ can also be ``truncated'' to a pair
$(\UresSL2, \HresSL2)$ of $2p^3$-dimensional algebras, where
$\HresSL2$ is a braided commutative Yetter--Drin\-feld
$\UresSL2$-module algebra.

\subsubsection{}
Interestingly, the $2p^3$-dimensional braided commutative
Yetter--Drin\-feld $\UresSL2$-module algebra $\HresSL2$ can be
described as
\begin{equation*}
  \HresSL2 \cong\Mat_p\bigl(\oC_{2p}[\lambda]\bigr),\qquad
  \oC_{2p}[\lambda]\equiv\oC[\lambda]/(\lambda^{2p}-1),
\end{equation*}
which adds a matrix flavor to our example.  In the matrix language,
the relevant structures are described as follows.  

First, the $\UresSL2$ action on matrices $X=(x_{ij})$ with
$\lambda$-dependent entries is given by
\begin{align}\label{K-mat}
  (K\leftii X)^{\vphantom{c}}_{ij}&
  {}=\q^{2(i-j)}
  \bigl(x^{\vphantom{c}}_{ij}\bigr|_{\lambda\to\q^{-1}\lambda}\bigr),
  \\[-6pt]
  \intertext{and}
  E\leftii X &= \smash[t]{\ffrac{1}{\q - \q^{-1}}
    \bigl(X\,\zM - \zM\,(K\leftii X)\bigr)},
  \\
  \label{F-mat}
  F\leftii X &= \ffrac{1}{\q - \q^{-1}}
  \bigl(\DzM\,X - (K^{-1}\leftii X)\,\DzM\bigr),
\end{align}
where
\begin{equation}
  \label{zDz-matrices}
  \zM={}
  \Matrixii{
    0&\hdotsfor{3}&0\\
    1&0&\hdotsfor{2}&0\\
    0&1&0&\dots&0\\
    \vdots&&\ddots&\ddots&\vdots\\
    0&\hdotsfor{2}&1&0
    },\qquad
  \DzM=(\q-\q^{-1})
  \Matrixii{
    0 & 1 & \hdotsfor{2} & \kern-12pt0\\
    0 & 0 & \;\q^{-1}\qint{2}\!\! & \hdots & \kern-12pt0\\
    \vdots&  & \ddots & \ddots & \kern-9pt\vdots \\
    0 & \hdotsfor{2}  &   0 & \q^{2-p}\qint{p-1}\\
    0 & \hdotsfor{3} & \kern-12pt 0
  }
\end{equation}
and we use the standard notation
\begin{equation*}
  \qint{n}=\ffrac{\q^n - \q^{-n}}{\q - \q^{-1}},\quad
  \qfac{n}=\qint{1}\qint{2}\dots\qint{n},\quad
  \qbin{m}{n}=\ffrac{\qfac{m}}{\qfac{m-n}\qfac{n}}.
\end{equation*}

Next, to describe the coaction
$\delta:\Mat_p(\oC_{2p}[\lambda])\to\UresSL2\tensor
\Mat_p(\oC_{2p}[\lambda])$, we first note that $\oC_{2p}[\lambda]$ is
the algebra of coinvariants, $\delta:\lambda\mapsto 1\tensor\lambda$.
It therefore remains to define $\delta$ on ``constant'' matrices
$\Mat_p(\oC)$.  But the full matrix algebra $\Mat_p(\oC)$ is
algebraically generated by the above $\zM$ and $\DzM$, and we have
\begin{equation}\label{coaction-M}
  \delta:\begin{aligned}
    \zM &\mapsto
    K^{-1}\tensor\zM
    -(\q - \q^{-1})
    E K^{-1}\tensor 1,
    \\
    \DzM &\mapsto
    K^{-1}\tensor \DzM
    +(\q - \q^{-1})
    F \tensor 1.
  \end{aligned}
\end{equation}
To summarize,
\begin{thm}
  With the above $\UresSL2$ action and coaction,
  $\Mat_p(\oC_{2p}[\lambda])$ and $\Mat_p(\oC)$ are braided
  commutative left--left Yetter--Drin\-feld $\UresSL2$-module
  algebras.
\end{thm}

\subsubsection{}By Theorem~4.1 in~\cite{[BM]}, as already noted
in~\bref{BM-1st}, we then have examples of bialgebroids:
\begin{equation*}
  \Mat_p(\oC_{2p}[\lambda])\Smash\UresSL2\ \ \text{and} \ \
  \Mat_p(\oC)\Smash\UresSL2
\end{equation*}
are Hopf algebroids over the respective algebras
$\Mat_p(\oC_{2p}[\lambda])$ and $\Mat_p(\oC)$; further details are
given in~\bref{algebroid}.

\subsection{Hopf algebras and logarithmic conformal field
  theory}\label{Rem:KLx} An additional source of interest in
$\UresSL2$ is its occurrence in a version of the Kazhdan--Lusztig
duality~\cite{[KLx]}, specifically, as the quantum group dual to a
class of \textit{logarithmic models} of conformal field
theory~\cite{[FGST],[FGST2],[S-q],[AM-3],[NT]}.

In the ``logarithmic'' Kazhdan--Lusztig duality, $\UresSL2$ appeared
in~\cite{[FGST],[FGST2]}; subsequently, it gradually transpired (with
the final picture having emerged from~\cite{[KS]}) that that was just
a continuation of a series of previous (re)discoveries of this quantum
group~\cite{[AGL],[Su],[X]} (also see~\cite{[Erd]}).  The ribbon and
(somewhat stretching the definition) factorizable structures of
$\UresSL2$ were worked out in~\cite{[FGST]}.

That $\UresSL2$ is Kazhdan--Lusztig-dual to logarithmic models of
conformal field theory\,---\,specifically, $\UresSL2$ at
$\q=e^{\frac{i\pi}{p}}$ is dual to the $(p,1)$ logarithmic
model~\cite{[Gaberdiel-K]}\,---\,means several things, in particular,
(i)~the $SL(2,\oZ)$ representation on the $\UresSL2$ center coincides
with the $SL(2,\oZ)$ representation generated from the characters of
the symmetry algebra of the logarithmic model~\cite{[FGST]}, the
so-called triplet $W(p)$
algebra~\cite{[Kausch],[Gaberdiel-K],[Gaberdiel-K-3],[FHST],[CF],
  [Adamovic-M-12]}, and (ii)~the $\UresSL2$ and $W(p)$ representation
categories are equivalent~\cite{[FGST2],[AM-3],[NT]}.

The ``Heisenberg counterpart'' of $\UresSL2$, its braided commutative
Yetter--Drin\-feld module algebra $\HresSL2$, is also likely to play a
role in the Kazhdan--Lusztig context~\cite{[S-U],[S-H]}, but this is a
subject of future work.

\section{$\HD(B^*)$ as a Yetter--Drin\-feld $\DD(B)$-module
  algebra}\label{sec:proofs}
We begin with simple facts about Yetter--Drinfeld module algebras,
concentrating in~\bref{sec:YD} on the construction of a braided
commutative Yetter--Drinfeld module algebra as a braided product
$X\Braid Y$ of two such algebras $X$ and $Y$.  In~\bref{want-use}, we
then specialize to $X=B^{*\mathrm{cop}}$ and $Y=B$, viewed as $\DD(B)$
module algebras under the heterotic action.  We verify that all the
necessary conditions are then satisfied, hence our conclusion
in~\bref{we-conclude}.  In~\bref{pseudo}, we give a ``pseudoadjoint''
interpretation of the heterotic action, and in~\bref{sec:multiple}
consider multiple ``alternating'' braided products.

\subsection{}\label{sec:YD}
The category of Yetter--Drin\-feld modules over a Hopf algebra with
bijective antipode is well known to be braided, with the braiding
$c^{\vphantom{1}}_{X,Y}:X\tensor Y\to Y\tensor X$ given by
\begin{equation*}
  c^{\vphantom{1}}_{X,Y}:x\tensor y\mapsto (x\mone\leftii y)\tensor x\zero.
\end{equation*}
The inverse is $c^{-1}_{X,Y}:y\tensor x\mapsto x\zero\tensor
S^{-1}(x\mone)\leftii y$.

We say that two Yetter--Drin\-feld modules $X$ and $Y$ are
\textit{braided symmetric} if
\begin{equation*}
  c^{\vphantom{1}}_{Y,X}=c_{X,Y}^{-1}
\end{equation*}
(note that both sides here are maps $Y\tensor X\to X\tensor Y$), that
is,
\begin{equation*}
  (y\mone\leftii x)\tensor y\zero
  = x\zero\tensor\bigl(S^{-1}(x\mone)\leftii y\bigr).
\end{equation*}
\begin{lemma}\label{lemma:locked}
  Let $X$ and $Y$ be braided symmetric Yetter--Drin\-feld modules, each
  of which is a braided commutative Yetter--Drin\-feld module algebra.
  Then their braided product $X\Braid Y$ is a braided commutative
  Yetter--Drin\-feld module algebra.
\end{lemma}

\subsubsection{Proof} 
Beyond the standard facts, we have to show the braided
commutativity,~i.e.,
\begin{equation}\label{to-show}
  \bigl((x\Braid y)\mone\leftii(v\Braid u)\bigr)(x\Braid y)\zero
  =(x\Braid y)(v\Braid u)
\end{equation}
for all $x,v\in X$ and $y,u\in Y$.  For this, we write the condition
$c^{\vphantom{1}}_{X,Y}=c^{-1}_{Y,X}$ as
\begin{equation*}
  (x\mone\leftii y)\tensor x\zero
  = y\zero\tensor\bigl(S^{-1}(y\mone)\leftii x\bigr)
\end{equation*}
and use this to establish an auxiliary identity,
\begin{equation}\label{locked}
  \bigl((x\mone\leftii y)\mone\leftii x\zero\bigr)
  \tensor(x\mone\leftii y)\zero
  \begin{aligned}[t]
    &=\bigl(y\zero{}\mone\leftii\bigl(S^{-1}(y\mone)\leftii x\bigr)
    \bigr)\tensor y\zero{}\zero
    \\
    &=\bigl(y''\mone S^{-1}(y'\mone)\leftii x\bigr)\tensor y\zero
    \\
    &=x\tensor y.
  \end{aligned}
\end{equation}
The left-hand side of~\eqref{to-show} can then be calculated as
\begin{align*}
  \bigl((x\Braid y)\mone\leftii(v\Braid u)\bigr)(x\Braid y)\zero
  \kern-140pt
  \\
    &=\bigl(x\mone y\mone\leftii(v\Braid u)\bigr)(x\zero\Braid y\zero)
    \\
    &=\bigl((x\mone' y\mone'\leftii v)\Braid
    (x\mone'' y\mone''\leftii u)\bigr)(x\zero\Braid y\zero)
    \\
    &=(x\mone' y\mone'\leftii v)
    \bigl((x\mone'' y\mone''\leftii u)\mone\leftii x\zero\bigr)
    \Braid
    (x\mone'' y\mone''\leftii u)\zero y\zero
    \\
    &=(x\mone y\mone'\leftii v)
    \bigl((x\zero{}\mone\leftii( y\mone''\leftii u))\mone\leftii
    x\zero{}\zero\bigr)
    \Braid
    (x\zero{}\mone\leftii(y\mone''\leftii u))\zero y\zero
    \\
    &=(x\mone y\mone'\leftii v)
    x\zero{}
    \Braid (y\mone''\leftii u)y\zero,
\end{align*}
just because of~\eqref{locked} in the last equality.  But the
right-hand side of~\eqref{to-show} is
\begin{align*}
  (x\Braid y)(v\Braid u)
  &=x (y\mone\leftii v) \Braid y\zero u
  \\[-3pt]
  &=(x\mone y\mone\leftii v)x\zero
  \Braid (y\zero{}\mone\leftii u) y\zero{}\zero
\end{align*}
because $X$ and $Y$ are both braided commutative.  The two expressions
coincide.

\subsubsection{Remark}\label{opp-braid}Because the braided symmetry
condition is symmetric with respect to the two modules, we also have
the braided symmetric Yetter--Drin\-feld module algebra $Y\Braid X$,
with the product
\begin{equation*}
  (y\Braid x)(u\Braid v)=y(x\mone\leftii u)\Braid x\zero v.
\end{equation*}
In addition to the multiplication inside $Y$ and inside $X$, this
formula expresses the relations $x u = (x\mone\leftii u)x\zero$
satisfied in $Y\Braid X$ by $x\in X$ and $u\in Y$.  Because
$c^{\vphantom{1}}_{X,Y}=c_{Y,X}^{-1}$, these are the same relations
$ux = (u\mone\leftii x)u\zero$ that we have in $X\Braid Y$.  Somewhat
more formally, the isomorphism
\begin{equation*}
  \phi: X\Braid Y\to Y\Braid X
\end{equation*}
is given by $\phi:x\Braid y\mapsto (x\mone\leftii y)\Braid x\zero$.
This is a module map by virtue of the Yetter--Drin\-feld condition, and
it is immediate to verify that $\delta(\phi(x\Braid
y))=(\id\tensor\phi)(\delta(x\Braid y))$.
That 
$\phi$ is an algebra map follows by calculating
\begin{align*}
  \phi(x\Braid y)\phi(v\Braid u)
  &=((x\mone\leftii y)\Braid x\zero)
  ((v\mone\leftii u)\Braid v\zero)
  \\
  &=(x\mone\leftii y)(x\zero{}\mone v\mone\leftii u)
  \Braid x\zero{}\zero v\zero
  \\
  &=(x'\mone\leftii y)(x''\mone v\mone\leftii u)
  \Braid x\zero v\zero
  \\
  &=x\mone\leftii\bigl(y(v\mone\leftii u)\bigr)
  \Braid x\zero v\zero\\
  \intertext{and}
  \phi((x\Braid y)(v\Braid u))
  &=\phi\bigl(x(y\mone\leftii v)\Braid y\zero u\bigr)
  \\
  &=
  (x\mone(y\mone\leftii v)\mone\leftii(y\zero u))
  \Braid x\zero(y\mone\leftii v)\zero
  \\
  &\stackrel{\checkmark}{=}
  x\mone\leftii(y\zero u)\zero\Braid
  x\zero\bigl(S^{-1}(y\zero{}\mone u\mone)\leftii(y\mone\leftii v)\bigr)
  \\
  &=
  x\mone \leftii(y\zero u\zero)\Braid
  x\zero\bigl(S^{-1}(y''\mone u\mone)y'\mone\leftii v\bigr)
  \\
  &=
  x\mone\leftii(y u\zero)\Braid
  x\zero\bigl(S^{-1}(u\mone)\leftii v\bigr)
  \\
  &\stackrel{\checkmark}{=}
  x\mone\leftii\bigl(y(v\mone\leftii u)\bigr)
  \Braid
  x\zero v\zero,
\end{align*}
where the braided symmetry condition was used in each of the
${}\stackrel{\checkmark}{=}{}$ equalities.

\subsection{}\label{want-use}We intend to use~\bref{lemma:locked} in
the case where $X=B^{*\mathrm{cop}}$ and $Y=B$.  This requires some
preparations.\enlargethispage{\baselineskip}

\begin{lemma}
  \label{lemma:2actions}
  For a Hopf algebra $B$ with bijective antipode, the formulas
  \begin{equation*}
    (\mu\tensor m)\leftii\beta=\mu''(m\leftreg\beta)\Sinv(\mu'),
    \qquad
    (\mu\tensor m)\leftii b
    =(m' b S(m''))\rightreg\Sinv(\mu)
  \end{equation*}
  make $B^{*\mathrm{cop}}$ and $B$ into left $\DD(B)$-module algebras.
\end{lemma}

\subsubsection{}\label{both}
This is known, e.g., from~\cite{[P]}, where both these actions are
discussed and references to the previous works are given.\pagebreak[3]
The $\DD(B)$ action on $B^*$ is obtained by restricting the
\textit{left} regular action of $\DD(B)$ on $\DD(B)^*\cong B\tensor
B^*$~\cite{[Lu-double]},
\begin{align*}
  (\mu\tensor m)\leftreg (b\tensor\beta)
  &= (\mu''\leftreg b)\tensor
  \mu'''(m\leftreg\beta)\Sinv(\mu'),
\end{align*}
to $1\tensor B^*$.  Similarly, the $\DD(B)$ action on $B$ is
obtained~\cite{[Z]} by restricting the \textit{right} regular action
of $\DD(B)$ on $\DD(B)^*\cong B\tensor B^*$ to $B\tensor\varepsilon$
and using the antipode to convert it into a left action.  The right
regular action of $\DD(B)$ on $\DD(B)^*$ is~\cite{[Lu-double],[P]}
\begin{equation*}
  (b\tensor\beta)\rightreg(\mu\tensor m)
  =S^{-1}(m''')(b\rightreg\mu)m'\tensor(\beta\rightreg m''),
\end{equation*}
where $\beta\rightreg m=\eval{\beta'}{m}\beta''$ is the right regular
action of $B$ on $B^*$.  Restricting to $B$ and replacing $\mu\tensor
m$ with $(S(m''')\leftreg\Sinv(\mu)\rightreg m')\tensor S(m'')$ then
gives the second formula in the lemma.

The following statement is obvious.
\begin{lemma}\label{lemma:2coactions}
  With the respective coactions
  \begin{equation*}
    \delta:\beta\mapsto(\beta''\tensor 1)\tensor\beta',
    \qquad
    \delta: b\mapsto (\varepsilon\tensor b')\tensor b'',
  \end{equation*}
  $B^{*\mathrm{cop}}$ and $B$ are $\DD(B)$-comodule algebras.
\end{lemma}

\begin{lemma}
  With the action and coaction in~\bref{lemma:2actions}
  and~\bref{lemma:2coactions}, both $B^{*\mathrm{cop}}$ and $B$ are
  Yetter--Drin\-feld module algebras.
\end{lemma}
It only remains to verify the Yetter--Drin\-feld condition in each case.
For $B^{*\mathrm{cop}}$, we calculate the left-hand side
of~\eqref{eq:YD} with $M=\mu\tensor m$ as
\begin{align*}
  ((\mu''\tensor m')\leftii\beta)\mone(\mu'\tensor m'')\tensor
  ((\mu''\tensor m')\leftii\beta)\zero\kern-190pt
  \\
  &=\bigl(\bigl(\mu'''(m'\leftreg\beta)\Sinv(\mu'')\bigr)''\mu'
  \tensor m''\bigr)
  \tensor
  (\mu'''(m'\leftreg\beta)\Sinv(\mu''))'
  \\
  &=\bigl(\mu\up{5}(m'\leftreg\beta'')\Sinv(\mu\up{2})\mu\up{1}
  \tensor m''\bigr)
  \tensor
  \mu\up{4}\beta'\Sinv(\mu\up{3})
  \\
  &=\bigl(\mu'''(m'\leftreg\beta'')
  \tensor m''\bigr)
  \tensor
  \mu''\beta'\Sinv(\mu'),
\end{align*}
but the right-hand side of~\eqref{eq:YD} is
\begin{align*}
  (\mu\tensor m)'\beta\mone\tensor
  \bigl((\mu\tensor m)''\leftii\beta\zero\bigr)\kern-60pt
  \\
  &=(\mu'''\tensor m')(\beta''\tensor 1)\tensor
  \bigl(\mu''(m''\leftreg\beta')\Sinv(\mu')\bigr)
  \\
  &=(\mu'''\tensor m')((\beta''\rightreg m'')\tensor 1)\tensor
  \mu''\beta'\Sinv(\mu')\\*[-3pt]
  &\phantom{{}={}}
  \text{\mbox{}\hfill(because
    $\beta''\tensor(m\leftreg\beta')=(\beta''\rightreg
    m)\tensor\beta'$)}
  \\[1pt]
  &=\bigr(\mu'''(m\up{1}\leftreg\beta''
  \rightreg m\up{4}S^{-1}(m\up{3}))\tensor m\up{2}\bigl)
  {}\tensor{}
  \mu''\beta'\Sinv(\mu'),
\end{align*}
which is the same.
For $B$, similarly, the left-hand side of~\eqref{eq:YD} is (assuming
the precedence $ab\rightreg\beta=(ab)\rightreg\beta$, and so on)
\begin{align*}
  ((\mu''\tensor m')\leftii b)\mone(\mu'\tensor m'')\tensor
  ((\mu''\tensor m')\leftii b)\zero\kern-220pt
  \\
  &=\bigl(\varepsilon\tensor
  \bigl((m' b S(m''))\rightreg\Sinv(\mu'')\bigr)'\bigr)
  (\mu'\tensor m''')\tensor
  \bigl((m' b S(m''))\rightreg\Sinv(\mu'')\bigr)''
  \\
  &=\bigl(\varepsilon\tensor
  \bigl((m' b S(m''))'\rightreg\Sinv(\mu'')\bigr)\bigr)
  (\mu'\tensor m''')\tensor
  (m' b S(m''))''
  \\*[-3pt]
  &\phantom{{}={}}
  \text{\mbox{}\hfill(because $\Delta(a\rightreg\mu)
    =(a'\rightreg\mu)\tensor a''$)}
  \\[1pt]
  &=\bigl(\mu''\tensor
  \bigl(\Sinv(\mu')\leftreg(m' b S(m''))'\bigr)m'''\bigr)\tensor
  (m' b S(m''))''\\*[-3pt]
  &\phantom{{}={}}
  \text{\mbox{}\hfill(using the $\DD(B)$-identity
    $\bigl(\varepsilon\tensor(b\rightreg\Sinv(\mu''))\bigr)(\mu'\tensor
    1) = \mu''\tensor(\Sinv(\mu')\leftreg b)$)}
  \\[1pt]
  &=\eval{\Sinv(\mu')}{m\up{2} b'' S(m\up{5})}
  \bigl(\mu''\tensor (m\up{1} b' S(m\up{6})) m\up{7}\bigr)
  \tensor
  m\up{3} b''' S(m\up{4})
  \\
  &=
  \bigl(\mu''\tensor m\up{1} b'\bigr)
  \tensor
  \bigl(m\up{2} b'' S(m\up{3})\rightreg\Sinv(\mu')\bigr)
  \\
  &=\bigl((\mu''\tensor m')(\varepsilon\tensor b')\bigr)\tensor
  \bigl((\mu'\tensor m'')\leftii b''\bigr)
  \\
  &=\bigl((\mu\tensor m)'b\mone\bigr)\tensor
  \bigl((\mu\tensor m)''\leftii b\zero\bigr),
\end{align*}
which is the right-hand side.

\begin{lemma}
  $B^{*\mathrm{cop}}$ and $B$ are braided commutative $\DD(B)$-module
  algebras.
\end{lemma}
This is entirely obvious once we note that when the $\DD(B)$ action on
$B^{*\mathrm{cop}}$ in~\bref{lemma:2actions} is restricted to the
action of $B^{*\mathrm{cop}}\tensor 1$, it becomes the adjoint action;
the same is true for the $\DD(B)$ action on $B$ restricted to the
action of $\varepsilon\tensor B$; therefore, for example,
$(a\mone\leftii b)a\zero=(a'\leftii b)a''=(a' b S(a''))a'''=a b$.

\begin{lemma}
  $B^{*\mathrm{cop}}$ and $B$ are braided symmetric.
\end{lemma}
We must show that
$c^{\vphantom{1}}_{B^{*\mathrm{cop}},B}=c_{B,B^{*\mathrm{cop}}}^{-1}$,
i.e.,
\begin{equation*}
  (b\mone\leftii\beta)\tensor b\zero
  =\beta\zero\tensor(S_{_{\DD}}^{-1}(\beta\mone)\leftii b).
\end{equation*}
The antipode here is that of $\DD(B)$, and therefore the right-hand
side evaluates as $%
\beta'\tensor(S^{*}(\beta'')\leftii b) =\beta'\tensor(b\rightreg
\Sinv(S^{*}(\beta''))) =\beta'\tensor(b\rightreg \beta'')$, which is
immediately seen to coincide with the left-hand side.

\subsection{}\label{we-conclude}
It now follows from~\bref{lemma:locked} that $B^{*\mathrm{cop}}\Braid
B$ is a braided commutative Yetter--Drin\-feld $\DD(B)$-module
algebra.  But the product in $B^{*\mathrm{cop}}\Braid B$ actually
evaluates as the product in~$\HD(B^*)$:
\begin{equation*}
  (\alpha\Braid a)(\beta\Braid b)=
  \alpha(a\mone\leftii\beta)\Braid a\zero b
  =\alpha((\varepsilon\tensor a')\leftii\beta)\Braid a''b
  =\alpha(a'\leftreg\beta)\Braid a''b.
\end{equation*}
We therefore conclude that \textit{with the $\DD(B)$ action and
  coaction in~\eqref{the-action} and~\eqref{delta}, $\HD(B^*)$ is a
  braided commutative Yetter--Drin\-feld $\DD(B)$-module algebra}.

\subsection{A ``pseudo-adjoint'' interpretation of the $\DD(B)$ action on
  $\HD(B^*)$}\label{pseudo} The action defined in~\eqref{the-action}
can be written in the ``pseudo-adjoint'' form
\begin{equation}\label{pseudo-ad}
  (\mu\tensor m)\leftii(\alpha\Smash a)=
  (\mu''\Smash m')\Star(\alpha\Smash a)
  \Star
  \pS(\mu'\tensor m''),
\end{equation}
where ${}\Star{}$ temporarily denotes the composition in~$\HD(B^*)$,
and
\begin{align*}
  \pS(\mu\tensor m)
  &=(\varepsilon\Smash S(m))\Star(\Sinv(\mu)\Smash 1)\\[-2pt]
  &=(S(m'')\leftreg\Sinv(\mu))\Smash S(m').
\end{align*}
The right-hand side of~\eqref{pseudo-ad} is to be compared with the
adjoint action of $\DD(B)$ on itself,
\begin{equation*}
  (\mu\tensor m)\adjoint(\nu\tensor n)=
  (\mu''\tensor m')\,(\nu\tensor n)\,S_{_{\DD(B)}}(\mu'\tensor m''),
\end{equation*}
where $S_{_{\DD(B)}}(\mu\tensor m)=(\varepsilon\tensor
S(m))(\Sinv(\mu)\tensor 1)
=(S(m''')\leftreg\Sinv(\mu)\rightreg m')\tensor S(m'')$.

\subsubsection{}To show~\eqref{pseudo-ad}, we calculate its right-hand
side as
\begin{align*}
  (\mu''\Smash m')\Star(\alpha\Smash a)\Star
  \bigl((S(m''')\leftreg\Sinv(\mu'))\Smash S(m'')\bigr)\kern-160pt
  \\[-1pt]
  &=\bigl(\mu''(m\up{1}\leftreg\alpha)\Smash m\up{2}a\bigr)
  \Star
  \bigl((S(m\up{4})\leftreg\Sinv(\mu'))\Smash S(m\up{3})\bigr)
  \\[-1pt]
  &=\mu''(m\up{1}\leftreg\alpha)
  \bigl(m\up{2} a' S(m\up{5})\leftreg\Sinv(\mu')\bigr)\Smash
  m\up{3} a'' S(m\up{4})
  \\[-1pt]
  &=\mu''(m'\leftreg\alpha)
  \bigl((m'' a S(m'''))'\leftreg\Sinv(\mu')\bigr)\Smash
  (m'' a S(m'''))''
  \\[-1pt]
  &=\mu'''(m'\leftreg\alpha)\Sinv(\mu'')
  \Smash
  (m'' a S(m''')\rightreg\Sinv(\mu'))
\end{align*}
(because
$(a'\leftreg\mu)\tensor a'' = \mu'\tensor(a\rightreg\mu'')$).

\subsubsection{}It may be interesting to see in more detail
\textit{why} the mock-adjoint action in~\eqref{pseudo-ad} is a
$\DD(B)$ action.  We recall from~\cite{[Lu-double]} that
Eq.~\eqref{Star-eta} holds for the product on~$\HD(B^*)$, with the
$\DD(B)$ product in the right-hand side and with the $2$-cocycle
$\eta:\DD(B)\tensor\DD(B)\to k$ given by
\begin{equation*}
  \eta(\mu\tensor m,\nu\tensor n)
  =\eval{\mu}{1}\eval{\nu}{m}\eval{\varepsilon}{n}.
\end{equation*}
Of course, $(M,A)\mapsto M\Star A$ is not a left action and
$(M,A)\mapsto A\Star \pS(M)$ is not a right action of $\DD(B)$;
instead, we have the associativity of the ${}\Star{}$ product,
$M\Star(A\Star N)=(M\Star A)\Star N$ for all $M,A,N\in B^*\tensor B$.
\ But the identity $\eta(M',N')\,\pS(N'')\Star\pS(M'') =\pS(M N)$
satisfied by Lu's cocycle~$\eta$ and the ``pseudo-antipode'' $\pS$
ensures that~\eqref{pseudo-ad} (i.e., \eqref{pseudo-ad-2}) is
nevertheless a $\DD(B)$ action.

From this perspective, furthermore, the $\DD(B)$ module algebra
property of $\HD(B^*)$ is ensured by another ``antipode-like''
property of~$\pS$, $\pS(M')\Star M''=\varepsilon(M)1$, $M\in\DD(B)$.
And the Yetter--Drin\-feld condition easily follows for the
``pseudo-adjoint'' action because $\delta\pS(M)=S(M'')\tensor\pS(M')$
(where $\delta$ is the same as~$\Delta_{_{\DD(B)}}$ and the right-hand
side is viewed as an element of $\DD(B)\tensor\HD(B^*)$) and, of
course, because $\HD(B^*)$ is a $\DD(B)$ comodule
algebra~\cite{[Lu-double]}.  We somewhat formalize this simple
argument as the following theorem (all of whose conditions hold for
Lu's cocycle).

\begin{thm}\label{thm:star}
  For a Hopf algebra $(H,\Delta,S,\varepsilon)$ with bijective
  antipode, let $\eta$ be a normal right
  $2$-cocycle~\cite{[Lu-double]}, i.e., a bilinear map $H\tensor H\to
  k$ such that
  \begin{equation*}
    \eta(f' g', h)\eta(f'', g'')=\eta(f, g' h')\eta(g'', h''),\quad
    \eta(1,h)=\eta(h,1)=\varepsilon(h)
  \end{equation*}
  for all $f,g,h\in H$, and let $H_{\Star}=(H,{}\Star{})$ denote the
  associative algebra with the product
  \begin{equation*}
    g\star h = g' h' \eta(g'', h'').
  \end{equation*}
  Let $\pS: H\to H$ be given by
  \begin{equation}\label{s-def}
    \pS(h)=\eta(h',h'')S(h''').
  \end{equation}
  If the conditions
  \begin{align}\label{etaX2}
    \eta(\pS(h'),h'')&=\varepsilon(h),
    \\
    \eta(h',\pS(h''))&=\varepsilon(h),
    \\
    \label{etaX5}
    \eta(g',h')
    \eta(\pS(h''),\pS(g''))&=\eta(g' h', g'' h'')
  \end{align}
  hold for all $g,h\in H$, then $H_{\Star}$ is a left--left
  Yetter--Drin\-feld $H$-module algebra under the left $H$-action
  \begin{equation}\label{thm-action}
    g\leftii h= g'\Star h\Star\pS(g'')
  \end{equation}
  and left coaction $\delta=\Delta$, viewed as a map $H_{\Star}\to
  H\tensor H_{\Star}$.  Moreover, $H_{\Star}$ is braided commutative.
\end{thm}
Conditions~\eqref{etaX2}--\eqref{etaX5} can be reformulated as
\begin{align}\label{thm-varepsilon}
  \pS(h')\Star h''
  &=\varepsilon(h)1,
  \\
  h' \Star\pS(h'')&=\varepsilon(h)1,
  \\
  \label{thm-cond}
  \eta(g',h')\pS(h'')\Star\pS(g'')&=\pS(g h).
\end{align}
Also, it follows from~\eqref{s-def} that
$\Delta(\pS(h)) = S(h'')\tensor\pS(h')$.
  
That~\eqref{thm-action} is an $H$ action immediately follows
from~\eqref{thm-cond}.  The module algebra property follows
from~\eqref{thm-varepsilon}.  The left coaction $\delta$ makes
$H_{\Star}$ into a comodule algebra for any right
cocycle~$\eta$~\cite{[Lu-double]}.  The Yetter--Drin\-feld axiom is then
verified as straightforwardly as for the true adjoint action:
\begin{align*}
  (h'\leftii g)\mone h''\tensor (h'\leftii g)\zero
  &= (h'\Star g\Star\pS(h''))' h''' \tensor (h'\Star g\Star\pS(h''))''
  \\
  &=h\up{1} g' S(h\up{4}) h\up{5}\tensor h\up{2}\Star g''\Star\pS(h\up{3})
  =h' g\mone\tensor (h''\leftii g\zero).
\end{align*}

The braided commutativity is also immediate:
\begin{align*}
  (h\mone\leftii g)\Star h\zero
  &=h'\mone\Star g\Star\pS(h''\mone)\Star h\zero
  =h'\Star g\Star\pS(h'')\Star h'''\\
  &=h'\Star g\Star 1\varepsilon(h'')=h\Star g.
\end{align*}

\subsection{Multiple braided products}\label{sec:multiple}
Further examples of Yetter--Drinfeld module algebras are produced by
extending the Heisenberg double~$\HD(B^*)$ to multiple ``alternating''
braided products.  We first return to the setting of~\bref{sec:YD}.

\subsubsection{}Multiple braided products $X_{1}\Braid\dots\Braid
X_{N}$ of Yetter--Drin\-feld $H$-module algebras $X_{i}$ are the
corresponding tensor products with the diagonal action and codiagonal
coaction of $H$, and with the relations\pagebreak[3]
\begin{equation}\label{i>j}
  x[i]\Braid y[j] = (x\mone\leftii y)[j]\Braid x\zero[i],
  \quad i>j,
\end{equation}
where $z[i]\in X_{i}$.  (The inverse relation is $x[i]\Braid y[j] =
y\zero[j]\Braid(S^{-1}(y\mone)\leftii x)[i]$, $i<j$.)  It readily
follows from the Yetter--Drin\-feld module algebra axioms for each of
the $X_i$ that $X_{1}\Braid\dots\Braid X_{N}$ is an associative
algebra and, in fact, a Yetter--Drin\-feld $H$-module algebra.  In
particular, it follows that
\begin{multline*}
  (x_1[i_1]\Braid\dots\Braid x_m[i_m])
  \Braid
  (y_1[j_1]\Braid\dots\Braid y_n[j_n])\\
  {}=
  \bigl((x_1{}\mone \dots x_m{}\mone)\leftii
  (y_1[j_1]\Braid\dots\Braid y_n[j_n])\bigr)\Braid
  \bigl(x_1{}\zero[i_1]\Braid\dots\Braid x_m{}\zero[i_m]\bigr)
\end{multline*}
whenever $i_a>j_b$ for all $a=1,\dots,m$ and $b=1,\dots,n$.

\subsubsection{``Alternating'' braided products}\label{XYXY}Next, let
$X$ and $Y$ be braided symmetric Yetter--Drin\-feld $H$-module
algebras, and consider the ``alternating'' products
\begin{equation*}
  X\Braid Y\Braid X\Braid Y\Braid\dots,
\end{equation*}
with an arbitrary number of factors (or a similar product with the
leftmost $Y$, or actually their inductive limits with respect to the
obvious embeddings).  We let $X[i]$ denote the $i$th copy of~$X$, and
similarly with~$Y[j]$.  For arbitrary $x[i]\in X[i]$ and $y[j]\in
Y[j]$, we then have relations~\eqref{i>j}, i.e.,
\begin{equation}\label{allij}
  x[2i+1]\Braid y[2j] = (x\mone\leftii y)[2j]\Braid x\zero[2i+1],
\end{equation}
for all $i\geq j$, but by the braided symmetry condition,
relations~\eqref{allij}\,---\,replicas of the relations between
elements of $X$ and elements of $Y$ in~\hbox{$X\Braid
  Y$}\,---\,\textit{hold for all $i$ and $j$}.  In the multiple
products, in addition, we also have the relations
\begin{equation}\label{xxyy}
  \begin{aligned}
    x[2i+1]\Braid v[2j+1]&=(x\mone\leftii v)[2j+1]\Braid x\zero[2i+1],
    \quad x,v\in X,\\[-2pt]
    y[2i]\Braid u[2j]&=(y\mone\leftii u)[2j]\Braid y\zero[2i],
    \quad y,u\in Y,
  \end{aligned}\quad i > j
\end{equation}
(which also hold for $i=j$ if $X$ and $Y$ are braided commutative.)

\subsubsection{Heisenberg $n$-tuples$/$chains}\label{sec:chains2}
Generalizing $\HD(B^*) \cong B^{*\mathrm{cop}}\Braid B\cong B \Braid
B^{*\mathrm{cop}}$, we have ``Heisenberg $n$-tuples$/$chains''---\,the
alternating products
\begin{align*}
  \HH_{2n} &=
  B^{*\mathrm{cop}}\Braid B\Braid B^{*\mathrm{cop}}\Braid B
  \Braid\dots\Braid B,\\[-2pt]
  \HH_{2n+1} &=
  B^{*\mathrm{cop}}\Braid B\Braid B^{*\mathrm{cop}}\Braid B
  \Braid \dots \Braid B\Braid B^{*\mathrm{cop}}.
\end{align*}
As we saw in~\bref{XYXY}, the following relations hold here:
\begin{alignat*}{2}
  b[2i]\,\beta[2j+1] &= (b'\leftreg\beta)[2j+1]\,b''[2i],
  \quad b\in B,\ \ \beta\in B^{*\mathrm{cop}},\ \ 
  \text{for all $i$ and $j$}\kern-200pt
  \\
  \intertext{(where $B^{*\mathrm{cop}}\to B^{*\mathrm{cop}}[2j+1]$ and
    $B\to B[2i]$ are the morphisms onto the respective factors, and we
    omit the ${}\Braid{}$ symbol for brevity), and}
  \alpha[2i+1]\,\beta[2j+1]&=(\alpha'''\beta\Sinv(\alpha''))[2j+1]\,
  \alpha'[2i+1],
  \quad&&\alpha,\beta\in B^{*\mathrm{cop}},\quad i\geq j,
  \\
  a[2i]\, b[2j]&=(a' b S(a''))[2j]\, a'''[2i],
  \quad&&a,b\in B, \quad i\geq j.
\end{alignat*}
The $\DD(B)$ action is diagonal and the coaction
is codiagonal, for example,
\begin{align*}
  \delta(\alpha\Braid a\Braid\beta\Braid b)
  &=\bigl((\alpha''\tensor 1)(\varepsilon\tensor
  a')(\beta''\tensor 1)(\varepsilon\tensor b')\bigr)
  \tensor
  \bigl(\alpha'\Braid a''\Braid\beta'\Braid b''\bigr)
  \\
  &=\bigl((\alpha''\tensor a')(\beta''\tensor b')\bigr)
  \tensor
  \bigl(\alpha'\Braid a''\Braid\beta'\Braid b''\bigr).
\end{align*}

The chains with the leftmost $B$ factor are defined entirely
similarly.  The obvious embeddings allow defining (one-sided or
two-sided) inductive limits of alternating chains.  All the chains are
Yetter--Drin\-feld module algebras, but none with $\geq3$ factors are
braided commutative in general.

\section{Yetter--Drinfeld module algebras and the associated Hopf
  algebroid for $\UresSL2$}\label{sec:sl2}
In this section, we construct Yetter--Drin\-feld module algebras for
$\UresSL2$ at the $2p$th root of unity for an integer $p\geq 2$
(see~\bref{sec:Uq}), and also consider the Hopf algebroid associated
with a braided commutative Yetter--Drinfeld module algebra in
accordance with the construction in~\cite{[BM]}.

$\UresSL2$ can be obtained as a subquotient of the Drin\-feld double
of a Taft Hopf algebra $B$~\cite{[FGST],[FGST2]} (a trick also used,
e.g., in~\cite{[Sch-Galois]} for a closely related quantum group).  On
the ``Heisenberg side,'' $\HD(B^*)$ similarly yields $\HresSL2$, a
$2p^3$-dimen\-sional braided commutative Yetter--Drin\-feld
$\UresSL2$-module algebra.  This is worked out
in~\bref{sec:HD}--\bref{sec:HD-red} below; in~\bref{sec:mat}, dropping
the coinvariants in $\HresSL2$, we obtain the algebra of $p\times p$
matrices, which is also a braided commutative Yetter--Drin\-feld
$\UresSL2$-module algebra.  In~\bref{algebroid}, we use the
Brzezi\'nski--Militaru theorem to construct the corresponding Hopf
algebroid.  Multiple alternating braided products are considered
in~\bref{sec:chains3}.

\subsection{$\DD(B)$ and $\HD(B^*)$ for the Taft Hopf algebra
  $B$}\label{sec:HD}
\subsubsection{The Taft Hopf algebra $B$}
Let 
\begin{equation*}
  B=\mathrm{Span}(E^m k^n),\quad
  0\leq m\leq p-1,\quad 0\leq n\leq 4p-1,
\end{equation*}
be the $4p^2$-dimensional Hopf algebra generated by~$E$ and~$k$ with
the relations
\begin{gather}\label{prod-B}
  k E =\q E k,\quad E^p=0,\quad k^{4p}=1,
\end{gather}
and with the comultiplication, counit, and antipode given by
\begin{gather}\label{coalgebra-B}
  \begin{gathered}
    \Delta(E)=1\otimes E+ E\otimes k^2,\quad
    \Delta(k)= k\otimes k,\quad
    \varepsilon(E)=0,\quad\varepsilon(k)=1,\\
    S(E)=- E k^{-2},\quad S(k)= k^{-1}.
  \end{gathered}
\end{gather}
We define $F,\varkappa\in B^*$ by
\begin{equation*}
  \eval{F}{E^{m} k^{n}}=\delta_{m,1}\ffrac{\q^{-n}}{\q-\q^{-1}},
  \qquad
  \eval{\varkappa}{E^{m} k^{n}}=\delta_{m,0}\q^{-n/2}.
\end{equation*}
Then~\cite{[FGST]}
\begin{equation*}
  B^*=\mathrm{Span}(F^a\varkappa^b),\quad
  0\leq a\leq p-1,\quad 0\leq b\leq 4p-1.
\end{equation*}

\subsubsection{The Drin\-feld double $\DD(B)$}
Direct calculation shows~\cite{[FGST]} that the Drin\-feld
double $\DD(B)$ is the Hopf algebra generated by $E$, $F$, $k$, and
$\varkappa$ with the relations given~by
\begin{itemize}
\item[i)] relations~\eqref{prod-B} in~$B$,
\item[ii)] the relations
    $\varkappa F=\q F\varkappa$, $F^p=0$,
    and $\varkappa^{4p}=1$
  in $B^*$, and
\item[iii)] the cross-relations
  \begin{gather*}
    k\varkappa=\varkappa k,\quad k F k^{-1}=\q^{-1} F,\quad
    \varkappa E\varkappa^{-1}=\q^{-1} E,\quad
    [E, F]=\mfrac{ k^2-\varkappa^2}{\q-\q^{-1}}.
  \end{gather*}
\end{itemize}
The Hopf-algebra structure
$(\Delta_{_{\DD}},\varepsilon_{_{\DD}},S_{_{\DD}})$ of $\DD(B)$ is
given by~\eqref{coalgebra-B} and
\begin{gather*}
  \Delta_{_{\DD}}(F)=\varkappa^2\tensor F+ F\tensor1,\quad
  \Delta_{_{\DD}}(\varkappa)=\varkappa\tensor\varkappa,\quad
  \varepsilon_{_{\DD}}(F)=0,\quad
  \varepsilon_{_{\DD}}(\varkappa)=1,
  \\
  S_{_{\DD}}(F)=-\varkappa^{-2} F,
  \quad S_{_{\DD}}(\varkappa)=\varkappa^{-1}.
\end{gather*}

\subsubsection{The Heisenberg double $\HD(B^*)$}\label{HD-relations}
For the above $B$, \ $\HD(B^*)$ is spanned by
\begin{equation}\label{HBstar-span}
  F^a\varkappa^b\Smash E^c k^d,\qquad a,c=0,\dots,p-1,\quad
  b,d\in\oZ/(4p\oZ),
\end{equation}
where $\varkappa^{4p}=1$, $k^{4p}=1$, $F^p=0$, and $E^p=0$.  A
convenient basis in $\HD(B^*)$ can be chosen as $(\varkappa, z,
\lambda, \Dz)$, where $\varkappa$ is understood as $\varkappa\Smash1$
and
\begin{align*}
  z &= -(\q-\q^{-1}) \varepsilon\Smash E k^{-2},\\
  \lambda &=\varkappa \Smash k,\\
  \Dz &= (\q-\q^{-1}) F\Smash 1.
\end{align*} 
The relations in $\HD(B^*)$ then become $\varkappa z = \q^{-1} z
\varkappa$, $\varkappa \lambda = \q^{\half} \lambda \varkappa$,
$\varkappa \Dz = \q \Dz \varkappa$, $\varkappa^{4p}=1$, and
\begin{equation}\label{Hsl2-rel}
  \begin{gathered}
    \lambda^{4p}=1,\qquad
    z^p=0,\qquad \Dz^p=0,
    \\
    \lambda z = z \lambda,\qquad   \lambda \Dz = \Dz \lambda,\\
    \Dz z = (\q-\q^{-1}) 1 +  \q^{-2} z\Dz.
  \end{gathered}
\end{equation}
Then the $\DD(B)$ action on $\HD(B^*)$ in~\eqref{the-action} becomes
$\varkappa\leftii\varkappa^n=\varkappa^n$,
\,$\varkappa\leftii\Dz^n=\q^{n}\Dz^n$,
\,$\varkappa\leftii\lambda^n=\q^{\frac{n}{2}}\lambda^n$,
\,$\varkappa\leftii z^n=\q^{-n}z^n$, and
\begin{equation}\label{D-action}
  \begin{alignedat}{3}
    E\leftii\varkappa &= 0,
    &
    k\leftii\varkappa^n &= \q^{-\frac{n}{2}}\varkappa,
    &
    F\leftii\varkappa^n&=-\q^{\frac{n}{2}}\qint{\fffrac{n}{2}}
    \dd\varkappa^n,
    \\
    E\leftii\lambda^n&=\q^{-\frac{n}{2}}\qint{\fffrac{n}{2}}\,\lambda^n\,z,
    &
    k\leftii\lambda^n&=\q^{-\frac{n}{2}}\lambda,
    &
    F\leftii\lambda^n&=-\q^{\frac{n}{2}}\qint{\fffrac{n}{2}}\,\lambda^n\,\Dz,
    \quad
    \\
    E\leftii z^n &=-\q^n \qint{n} z^{n+1},
    &
    k\leftii z^n &= \q^{n}\,z^n,
    &
    F\leftii z^n &= \qint{n} \q^{1-n}\,z^{n-1},
    \\
    E\leftii \Dz^n &= \q^{1-n}\qint{n}\Dz^{n-1},\quad
    &
    k\leftii \Dz^n &= \q^{-n} \Dz^n,\quad
    &
    F\leftii \Dz^n &= -\q^n \qint{n} \Dz^{n+1}.
  \end{alignedat}
\end{equation}

\subsection{The $(\UresSL2,\HresSL2)$ pair}\label{sec:HD-red}
\subsubsection{From $\DD(B)$ to $\UresSL2$}
The ``truncation'' whereby $\DD(B)$\pagebreak[3] yields $\UresSL2$
consists of two steps~\cite{[FGST]}: first, taking the quotient
\begin{gather}\label{quotient}
  \bDDB=\DD(B)/(\varkappa k - 1)
\end{gather}
by the Hopf ideal generated by the central element $\varkappa\tensor k
- \varepsilon\tensor 1$ and, second, identifying $\UresSL2$ as the
subalgebra in $\bDDB$ spanned by $F^{\ell} E^{m} k^{2n}$ with
$\ell,m=0,\dots,p-1$ and $n=0,\dots,2p-1$.  It then follows that
$\UresSL2$ is a Hopf algebra\,---\,the one described
in~\bref{sec:Uq}, where $K=k^2$.

The category of finite-dimensional $\UresSL2$ representations is not
braided~\cite{[KS]}.

\subsubsection{From $\HD(B^*)$ to $\HresSL2$}\label{sec:Hq}
In $\HD(B^*)$, dually to the two steps just mentioned, we take a
subalgebra and then a quotient~\cite{[S-H]}.  In the basis chosen
above, the subalgebra (which is also a $\UresSL2$ submodule) is the
one generated by $z$, $\Dz$, and $\lambda$.  Its quotient by
$\lambda^{2p}=1$ gives the $2p^3$-dimensional algebra
\begin{equation*}
  \HresSL2 = \oC[z,\Dz,\lambda]\!\bigm/\!(\eqref{Hsl2-rel}\text{ and }
  (\lambda^{2p}-1)).
\end{equation*}

As an associative algebra,
\begin{equation*}
  \HresSL2 = \Czd\tensor(\oC[\lambda]/(\lambda^{2p}-1)),
\end{equation*}
with the $p^2$-dimensional algebra
\begin{equation}\label{Czd-def}
  \Czd=\oC[z,\Dz]/(z^p, \ \Dz^p, \ \Dz z - (\q-\q^{-1}) - \q^{-2} z\Dz).
\end{equation}
The $\UresSL2$ action on $\HresSL2$ is given by the last three lines
in~\eqref{D-action}, with the central column rewritten for $K=k^2$.
The coaction $\delta:\HresSL2\to\UresSL2\tensor\HresSL2$ follows
from~\eqref{delta} as
\begin{align*}
  \lambda&\mapsto1\tensor\lambda,\\
  z^m &\mapsto \sum_{s=0}^m(-1)^s\q^{s(1-m)}(\q - \q^{-1})^s
  \qbin{m}{s}\,
  E^s K^{-m}\tensor z^{m-s},\\
  \Dz^m &\mapsto \sum_{s=0}^m\q^{s(m-s)}(\q - \q^{-1})^s
  \qbin{m}{s}\,
  F^s K^{s-m}\tensor \Dz^{m-s}.
\end{align*}

\subsubsection{} \textit{With the $\UresSL2$ action and coaction given
  above, $\HresSL2$ is a braided commutative Yetter--Drin\-feld
  $\,\UresSL2$-module algebra}.

\subsection{Matrix braided commutative Yetter--Drin\-feld module
  algebras}\label{sec:mat}
It follows that $\Czd$ in~\eqref{Czd-def}\,---\,the algebra of
``quantum differential operators on a line''---\,is also a braided
commutative Yetter--Drin\-feld $\,\UresSL2$-module algebra.  It is in
fact the full matrix algebra~\cite{[S-U]},
\begin{equation}\label{iso}
  \Czd\cong\Mat_p(\oC).
\end{equation}

\subsubsection{}
That $\Czd$ is (semisimple and) isomorphic to $\Mat_p(\oC)$ already
follows from a more general picture elegantly developed
in~\cite{[FIK]}, where ``para-Grassmann'' algebras of the
form $\oC[z,\Dz]/(z^{p},\Dz^{p})$ with \textit{various} additional
relations on the $z^i \Dz^j$ were studied.  The relations between our
$z$ and $\Dz$,
\begin{gather*}
  \Dz^m\,z^n=\sum_{i\geq0}\q^{-(2 m - i) n + i m - \frac{i(i-1)}{2}}
  \qbin{m}{i} \qbin{n}{i} \qfac{i} \left(\q - \q^{-1}\right)^i
  z^{n-i}\Dz^{m-i}
\end{gather*}
(where the range of $i$ is bounded above by $\min(m,n)$ because of the
$\q$-binomial coefficients), are nondegenerate in terms of the
classification in~\cite{[FIK]}, hence the isomorphism with the full
matrix algebra.

We describe~\eqref{iso} as an isomorphism \textit{of $\UresSL2$ module
  comodule algebras}.  The generators $z$ and $\Dz$ have the
respective matrix representations $\zM$ and $\DzM$
in~\eqref{zDz-matrices} (where we do not reduce the expressions using
that $\q^{p}=-1$ and $\qint{p-i}=\qint{i}$ to highlight a pattern).
Coaction~\eqref{coaction-M} is then just the $m=1$ case of the
formulas in~\bref{sec:Hq}, and it is not difficult to see that the
last three lines in~\eqref{D-action} yield
formulas~\eqref{K-mat}--\eqref{F-mat}\,---\,so far, with no effect of
the rescaling of $\lambda$~in~\eqref{K-mat}.

\subsubsection{}Once $\Czd$ is thus identified with $\Mat_p(\oC)$, we
can write
\begin{equation*}
  \HresSL2=\Mat_p\bigl(\oC_{2p}[\lambda]\bigr)
\end{equation*}
(where we recall that
$\oC_{2p}[\lambda]=\oC[\lambda]/(\lambda^{2p}-1)$), and it is
immediate to see from~\eqref{D-action} that $\lambda$ entering the
matrix entries rescales under the $\UresSL2$ action as indicated
in~\eqref{K-mat}.  This establishes
formulas~\eqref{K-mat}--\eqref{F-mat}.

\subsubsection{}For example, for $p=3$, choosing
$x_{ij}=\lambda^{n_{ij}} y_{ij}$ with $\lambda$-independent $y_{ij}$,
we have
\begin{multline*}
  F\leftii\Matrixi{
    \lambda ^{n_{11}} y_{11} & \lambda
   ^{n_{12}} y_{12} & \lambda ^{n_{13}}
   y_{13} \\
 \lambda ^{n_{21}} y_{21} & \lambda
   ^{n_{22}} y_{22} & \lambda ^{n_{23}}
   y_{23} \\
 \lambda ^{n_{31}} y_{31} & \lambda
   ^{n_{32}} y_{32} & \lambda ^{n_{33}}
   y_{33}
 }
  \\[-2pt]
  ={}
  \Matrixi{
    \lambda^{n_{21}} y_{21} & \lambda^{n_{22}} y_{22} -
    \q^{n_{11}}
    \lambda^{n_{11}} y_{11} &
    \lambda^{n_{23}} y_{23} + \q^{n_{12}-2} \lambda^{n_{12}} y_{12}\\[2pt]
    \q^{-1}\lambda^{n_{31}} y_{31} &
    \q^{-1} \lambda^{n_{32}} y_{32}
    - \q^{n_{21}-2} \lambda^{n_{21}} y_{21} &
    \q^{-1} \lambda^{n_{33}} y_{33}
    + \q^{n_{22}-4} \lambda^{n_{22}} y_{22}
    \\[2pt]
    0 & -\q^{n_{31}+2} \lambda^{n_{31}} y_{31} &
    \q^{n_{32}} \lambda^{n_{32}} y_{32}
  }.
\end{multline*}

\subsubsection{}As an example of the coaction in matrix form,
$\delta:\Mat_p(\oC)\to\UresSL2\tensor\Mat_p(\oC)$, we give the only
typographically manageable case, that of~$p=2$.  Writing elements of
$\UresSL2\tensor\Mat_2(\oC)$ as matrices with $\UresSL2$-valued
entries, we have
\begin{multline*}
  \delta X
  \\[-2pt]
  {}=\mbox{\footnotesize$\begin{pmatrix}
      (1-2 i
      E F K^3)
      x_{11}+F x_{12}-2 i
      E K^3 x_{21}+2 i
      E F K^3
      x_{22}\kern-45pt &\ 2 i E K^2
      x_{11}+K^3 x_{12}-2 i
      E K^2 x_{22}
      \\[4pt]
      F K^3 x_{11}+K^3
      x_{21}-F K^3 x_{22}\qquad &
      \kern-35pt(1-K^2-2 i
      E F K^3)
      x_{11}+F x_{12}-2 i
      E K^3 x_{21}+(K^2+2 i
      E F K^3)
      x_{22}
    \end{pmatrix}$}.
\end{multline*}

\subsubsection{}It would be interesting to find a direct
\textit{matrix} derivation of the Yetter--Drin\-feld axiom for
$\Mat_p(\oC)$ and the braided commutativity property
\begin{equation*}
  (X\mone\leftii Y)X\zero=XY,\qquad X,Y\in\Mat_p(\oC).
\end{equation*}
We illustrate the structure occurring in the left-hand
side here before the matrix multiplication, with the known result, is
evaluated (again, necessarily restricting ourself to $p=2$):%
\begin{multline*}\renewcommand{\arraycolsep}{4pt}%
  (X\mone\leftii Y)\tensor X\zero
  =
  \Matrixii{
    y_{11} & -y_{12} \\
    -y_{21} & y_{22}
  }
  \tensor
  \Matrixii{
    0 & x_{12} \\
    x_{21} & 0
  }
  +
  \Matrixii{
    y_{11} & y_{12} \\
    y_{21} & y_{22}
  }
  \tensor
  \Matrixii{
    x_{11} & 0 \\
    0 & x_{22}
  }
  \\
  {}+
  \Matrixii{
    -\frac{i}{2} y_{12} & 0 \\
    \frac{i}{2} (y_{11} - y_{22}) & -\frac{i}{2} y_{12}
  }
  \tensor
  \Matrixii{
    0 & 2 i (x_{11} - x_{22}) \\
    0 & 0
  }
  +
  \Matrixii{
    \frac{i}{2} y_{12} & 0 \\
    \frac{i}{2} (y_{11} - y_{22}) & \frac{i}{2} y_{12}
  }
  \tensor
  \Matrixii{
    -2 i x_{21} & 0 \\
    0 & -2 i x_{21}
  }
  \\
  {}+
  \Matrixii{
    y_{21} & y_{11}-y_{22} \\
    0 & y_{21}
  }
  \tensor
  \Matrixii{
    0 & 0 \\
    x_{22}-x_{11} & 0
  }
  +
  \Matrixii{
    y_{21} & y_{22}-y_{11} \\
    0 & y_{21}
  }
  \tensor
  \Matrixii{
    x_{12} & 0 \\
    0 & x_{12}
  }
  \\
  {}+
  \Matrixii{
    \frac{i}{2} (y_{22} - y_{11}) & 0 \\
    0 & \frac{i}{2} (y_{22} - y_{11})
  }
  \tensor
  \Matrixii{
    2 i (x_{11} - x_{22}) & 0 \\
    0 & 2 i (x_{11} - x_{22})
  }.
\end{multline*}

\subsection{Hopf algebroid with the $\Mat_p(\oC)$
  base}\label{algebroid}Theorem~4.1 in~\cite{[BM]} nicely
reinterprets the structure of a braided commutative Yetter--Drin\-feld
$H$-module algebra $A$ as a bialgebroid structure on $A\Smash H$.  (We
refer the reader to~\cite{[BM]} for a comprehensive discussion of
$($Hopf$|$bi$)$algebroids, also in relation to Lu's
bialgebroids~\cite{[Lu-alg]}, Xu's bialgebroids with an
anchor~\cite{[Xu]}, and Takeuchi's $\times_{A}$-bialgebras~\cite{[T]},
as well as for references to other related works.)

The examples of Hopf algebroids $A\Smash H$ with our braided
commutative Yetter--Drin\-feld module algebras
$A=\Mat_p(\oC_{2p}[\lambda])$ or $A=\Mat_p(\oC)$ may be of some
interest because of the explicit matrix structure of the base
algebra~$A$.  Below, we follow~\cite{[BM]}, adapting the formulas
there to a \textit{left} comodule algebra by duly inserting the
antipodes.  To somewhat simplify the notation, we discuss the
``$\lambda$-independent'' example, i.e., the Hopf algebroid structure
of $\XX=\Mat_p(\oC)\Smash\UresSL2$; reintroducing $\oC_{2p}[\lambda]$
on the matrix side is left to the reader.

As a vector space, $\XX\cong\Mat_p(\UresSL2)$, matrices with
$\UresSL2$-valued entries; we can therefore write $1\Smash h = \one h$
\ ($h\in\UresSL2$), where $\one$ is the unit $p\times p$ matrix; with
a slight abuse of notation, similarly, $X\Smash 1=X$, understood as a
``constant'' $p\times p$ matrix.  An arbitrary element of $\XX$ can be
written as $\sum_{i,j=1}^p \E_{ij}h_{ij}$, where the $\E_{ij}$ are the
standard elementary matrices and $h_{ij}\in\UresSL2$.  The
smash-product composition is then given by
\begin{equation*}
  \Bigl(\sum_{i,j=1}^p \E_{ij}h_{ij}\Bigr)
  \Bigl(\sum_{m,n=1}^p \E_{mn}g_{mn}\Bigr)
  =\sum_{i,j=1}^p\sum_{m,n=1}^p \E_{ij}
  (h_{ij}'\leftii\E_{mn})\, h_{ij}'' g_{mn},
  \quad h_{ij},g_{mn}\in\UresSL2,
\end{equation*}
with the left action $\leftii$ to be evaluated in accordance
with~\eqref{K-mat}--\eqref{F-mat}.  We write
$\XX=\Mat_p(\UresSL2)_{\Smashfn}$, with the subscript reminding of the
smash-product composition in this algebra (which is highly nonstandard
from the matrix standpoint).

The relevant structures
\begin{align*}
  \epsilon&:\XX\to\Mat_p(\oC),
  \\
  s,t&:\Mat_p(\oC)\to\XX,
  \\
  \Delta&:\XX\to\XX\tensorA\XX,
  \\
  \tau&:\XX\to\XX
\end{align*}
(the counit, the source and target maps, the coproduct, and the
antipode) are as follows.

The counit $\epsilon:\XX\to\Mat_p(\oC)$ acts componentwise,
\begin{equation*}
  \epsilon\Bigl(\sum_{i,j=1}^p \E_{ij}h_{ij}\Bigr)
  =\sum_{i,j=1}^p \E_{ij}\varepsilon(h_{ij}).
\end{equation*}

The source map $s:\Mat_p(\oC)\to\XX$ is the identical map onto
constant matrices.  The target map
$t:\Mat_p(\oC)\to\XX=\Mat_p(\oC)\Smash\UresSL2$ is given by
\begin{equation*}
  t(X)=X\zero\Smash S^{-1}(X\mone),
\end{equation*}
where $\delta(X)=X\mone\tensor X\zero\in\UresSL2\tensor\Mat_p(\oC)$ is
the coaction defined in~\eqref{coaction-M}.  It then follows that
$\zM$ and $\DzM$ in~\eqref{zDz-matrices} map under~$t$ into the
following two-diagonal matrices with $\UresSL2$-valued entries:
\begin{align}\label{t(Z)}
  t(\zM)&=
  \Matrixi{
    (\q-\q^{-1})E\kern-12pt&0\\
    K&\kern-12pt(\q-\q^{-1})E&0\\
    \vdots&\kern-40pt\ddots&\kern-10pt\ddots&\ddots
    \\
    0&\kern-30pt\hdots&\kern-20pt K\qquad&\kern-20pt(\q-\q^{-1})E&
    \kern-20pt 0\\
    0&\hdotsfor{2}&\kern-30pt K&\kern-10pt(\q-\q^{-1})E
  },
  \\
  \label{t(D)}
  t(\DzM)&=
  (\q-\q^{-1})
  \Matrixi{
    -FK&K&\\
    0&-FK&\q^{-1}\qint{2}K\\
    \vdots&\kern-30pt\ddots&\kern-30pt\ddots&\kern-10pt\ddots\\
    0&\hdots&\kern-40pt 0&\kern-30pt{-FK}\quad&\kern-10pt\q^{2-p}\qint{p-1}K\\
    0&\hdotsfor{1}&&\kern-40pt 0&\kern-8pt{-FK}\qquad
  }.
\end{align}
For any complex matrix $Y=\sum_{m,n}y_{mn}\zM^m\DzM^n$, we
use~\eqref{t(Z)} and~\eqref{t(D)} to calculate
$t(Y)=\sum_{m,n}y_{mn}t(\DzM)^n t(\zM)^m\in
\XX$ (evidently, with the smash-product multiplication understood).
Furthermore, elementary calculation using the braided commutativity
shows that 
\begin{equation*}
  t(Y)(X\Smash h) = X\cdot t(Y)\cdot h,\quad X,Y\in\Mat_p(\oC),
\end{equation*}
where, abusing the notation, the dot denotes both \textit{matrix
  product and the product in $\UresSL2$}, with articulately no
``smash'' effects because multiplication with a constant matrix is on
the left and with a $\UresSL2$ element on the right.

The coproduct $\DDelta:\XX\to\XX\tensorA\XX$ is (co)componentwise,
\begin{equation*}
  \DDelta\Bigl(\sum_{i,j=1}^p \E_{ij}h_{ij}\Bigr)
  =\sum_{i,j=1}^p \E_{ij}h'_{ij}\tensorA \one h''_{ij},
  \quad h_{ij}\in\UresSL2.
\end{equation*}
The $\tensorA$ product is here defined with respect to the right
action of $\Mat_p(\oC)$ on $\XX$ via $(X\Smash h).Y = t(Y)(X\Smash h)
$ and the left action via $Y.(X\Smash
h)=s(Y)(X\Smash h)
$, and hence
\begin{equation*}
  (X\cdot t(A)\cdot h)\tensorA(Y\Smash g)
  =(X\Smash h)\tensorA(A Y\Smash g)
\end{equation*}
holds for all $X,A,Y\in\Mat_p(\oC)$ and $g,h\in\UresSL2$ (where the
first factor in the left-hand side, again, involves matrix and
$\UresSL2$ products on the different sides of
$t(A)$).

The antipode $\tau:\XX\to\XX$ is given by another simple adaptation of
a formula in~\cite{[BM]}:
\begin{equation*}
  \tau(X\Smash h) = (1\Smash S(h))
  \bigl((S(X\mone'')\leftii X\zero)\Smash S(X\mone')\bigr),
  \quad X\in\Mat_p(\oC),\ \ h\in\UresSL2,
\end{equation*}
with the product in the right-hand side to be taken
in~$\XX$.\footnote{And the section $\gamma$ of the natural projection
  $\XX\tensor\XX\to\XX\tensorA\XX$, required in the definition of a
  Hopf algebroid~\cite{[BM]} to satisfy the condition
  $m\circ(\id\tensor\tau)\circ\gamma\circ\Delta = s\circ\epsilon$, is
  given by $\gamma:(X\Smash h)\tensorA(Y\Smash g) \mapsto
  (X\cdot t(Y)\cdot h)\tensor(\one g)$.}  On $1\Smash\UresSL2$, this
is just the $\UresSL2$ antipode, and on $\Mat_p(\oC)$, $\tau(X) =
(1\Smash S(X\mone))(X\zero\Smash 1)$; a simple calculation then shows
that
\begin{align*}
  \tau(\zM)&=\q^2 t(\zM),
  \\
  \tau(\DzM)&=\q^{-2} t(\DzM).
\end{align*}
Being an anti-algebra map, again, this extends to all of
$\Mat_p(\oC)\ni\sum_{m,n}y_{mn}\zM^m\DzM^n$.

Some of the Hopf algebroid properties (see~\cite[Defnition~2.2]{[BM]}
for a nicely refined list of axioms), e.g., $\tau(t(X))=s(X)$ and
$t(X)s(Y)=s(Y)t(X)$, are evident for $\XX=\Mat_p(\oC)\Smash\UresSL2$
described in matrix form; with others, it is not entirely obvious how
far one can proceed with verifying them in a purely \textit{matrix}
language, i.e., not following~\cite{[BM]} in resorting to the
Yetter--Drin\-feld module algebra properties; so much more interesting
is the fact that $\Mat_p(\UresSL2)_{\Smashfn}$ at $\q=e^{i\pi/p}$\,
\textit{is}\, a Hopf algebroid over~$\Mat_p(\oC)$.

As already noted, it is entirely straightforward to extend the above
formulas to describe $\Mat_p(\oC_{2p}[\lambda])\Smash\UresSL2
=\Mat_p(\UresSL2\tensor\oC_{2p}[\lambda])_{\Smashfn}$ as a Hopf
algebroid over
$\Mat_p(\oC_{2p}[\lambda])\equiv\Mat_p(\oC[\lambda]/(\lambda^{2p}-1))$.

\subsection{Heisenberg ``chains.''}\label{sec:chains3}
The Heisenberg $n$-tuples$/$chains defined in~\bref{sec:chains2} can
also be ``truncated''\pagebreak[3] similarly to how we passed from
$\HD(B^*)$ to $\HresSL2$.  An additional possibility here is to drop
the coinvariant $\lambda$ altogether, which leaves us with the
``\textit{truly Heisenberg}'' Yetter--Drin\-feld $\UresSL2$-module
algebras
\begin{align*}
  \HHH_2 &= \oC_{\q}^{*p}[\Dz_1]\Braid\oC_{\q}^p[z_2]
  = \oC_{\q}[z_2,\Dz_1],\\
  \HHH_{2n} &= \oC_{\q}^{*p}[\Dz_1]\Braid\oC_{\q}^p[z_2]\Braid\dots
  \Braid\oC_{\q}^{*p}[\Dz_{2n-1}]\Braid\oC_{\q}^p[z_{2n}],\\
  \HHH_{2n+1} &= \oC_{\q}^{*p}[\Dz_1]\Braid\oC_{\q}^p[z_2]\Braid\dots
  \Braid\oC_{\q}^{*p}[\Dz_{2n-1}]\Braid\oC_{\q}^p[z_{2n}]
  \Braid\oC_{\q}^{*p}[\Dz_{2n+1}]
\end{align*}
(or their infinite versions), where
$\oC_{\q}^{*p}[\Dz]=\oC[\Dz]/\Dz^p$ and $\oC_{\q}^{p}[z]=\oC[z]/z^p$
with the braiding inherited from~\bref{sec:chains2}, which amounts to
using the relations
\begin{align*}
  &\Dz_i\, z_{j} = \q - \q^{-1} + \q^{-2} z_{j}\,\Dz_i
  \\[-4pt]
  \intertext{for \textit{all} (odd) $i$ and (even) $j$, and}
  &\begin{aligned}
    z_i\,z_j &= \smash[t]{\q^{-2} z_j\,z_i + (1 - \q^{-2}) z_j^2},\\
    \Dz_i\,\Dz_j &= \q^2 \Dz_j\,\Dz_i + (1 - \q^2)\Dz_j^2,
  \end{aligned}\quad i\geq j
\end{align*}
(and $z_i^p=$ and $\Dz_i^p=0$; our relations may be interestingly
compared with those in para-Grassmann algebras studied
in~\cite{[FIK1]}).

\subsubsection*{Acknowledgments} I am grateful to J.~Fjelstad,
J.~Fuchs, and A.~Isaev for the useful comments.  This work was
supported in part by the RFBR grant 08-01-00737, the RFBR--CNRS grant
09-01-93105, and the grant LSS-1615.2008.2.

\parindent=0pt

\end{document}